\documentclass[12pt]{amsart}
\usepackage{epsfig}
\usepackage{amssymb,latexsym,bbm,graphicx,epsfig,epic,eepic,oldgerm}
\usepackage{amscd}
\usepackage[mathscr]{eucal}
\usepackage{amssymb}
\usepackage{amsxtra}
\usepackage{amsmath}
\usepackage[all]{xy}

\newtheorem{thm}{Theorem}
\newtheorem*{rem}{Remark}

\def\ind{\operatorname{ind}}
\def\Ker{\operatorname{Ker}}

\begin{document}

\title[The monotone Audin conjecture]{The Maslov class of Lagrangian tori and quantum products in Floer cohomology}
\author{Lev Buhovsky}
\address{Lev Buhovsky, School of Mathematical Sciences, Tel-Aviv
  University, Ramat-Aviv, Tel-Aviv 69978, Israel}
\email{levbuh@post.tau.ac.il}
\thanks{This is part of my PhD thesis, being carried out
under the guidance of Prof. P.~Biran, at Tel-Aviv University.}

\thanks{The author was partially supported by the ISRAEL SCIENCE FOUNDATION (grant No. 1227/06 *)}

\maketitle

\begin{abstract}

 We use Floer cohomology to prove the monotone version of a conjecture
of Audin: the minimal Maslov number of a monotone Lagrangian torus in
$\mathbb{R}^{2n}$ is $2$. Our approach is based on the study of the
quantum cup product on Floer cohomology and in particular the
behaviour of Oh's spectral sequence with respect to this product.  As
further applications we prove existence of holomorphic disks with
boundaries on Lagrangians as well as new results on Lagrangian
intersections.

\end{abstract}

\section{Introduction and main results} \label{S:intro}

 Let $(M, \omega)$ be a tame symplectic manifold (see ~\cite{A-L-P}, also Section ~\ref{S:basic notions}).
The class of tame symplectic manifolds includes compact manifolds,
Stein manifolds, and more generally, manifolds which are
symplectically convex at infinity, as well as products of all the
above. Let $L \subset (M, \omega)$ be a closed Lagrangian
submanifold. Throughout the paper, by a closed manifold we mean
compact manifold without boundary, and all Lagrangian manifolds
are supposed to be compact and without boundary. One of
fundamental problems in symplectic topology is to find
restrictions on the topology of $ L $, in particular on the Maslov
class $\mu_{L}:\pi_{2}(M,L) \rightarrow \mathbb{Z}$. See
Section~\ref{S:basic notions} for the definition, and for other
basic notions from symplectic topology. Below we will be concerned
with monotone Lagrangian submanifolds. For such a Lagrangian $ L
\subset (M, \omega) $ denote by $ N_L \in \mathbb{Z}_+ $ the
minimal Maslov number, i.e. the generator of the image of $
\mu_{L} $.

 Our first result deals with Lagrangian submanifolds of $ \mathbb{R}^{2n} $.
Here we endow $ \mathbb{R}^{2n} $ with its standard symplectic
structure $ \omega_{std} = dx_1 \wedge dy_1 + dx_2 \wedge dy_2 +
... + dx_n \wedge dy_n $, where $ (x_1,y_1,x_2,y_2,...,x_n,y_n) $  are standard coordinates on $ \mathbb{R}^{2n} $.

\begin{thm} \label{T:Audin conj}
   Given a monotone Lagrangian embedding
   of the $ n $-torus $ L = \mathbb{T}^{n} \hookrightarrow (
   \mathbb{R}^{2n}, \omega_{std} ) $, its minimal Maslov number must be
   $ N_{L} = 2 $.
\end{thm}

 Let us remark, that the minimal Maslov of such an embedding
must be even, due to the orientability of the torus, and it is a
non-negative integer. The possibility of $ N_L = 0 $ cannot occur,
since then, by monotonicity of $ L $, the area class of $ L $ must
vanish, and this contradicts the famous result of Gromov
~\cite{Gr}, that in our case guarantees the existence of a
pseudo-holomorphic disc with a boundary on $ L $, which must have
a positive symplectic area.

\begin{thm} \label{T:Maslov 2}

   Let $ L = \mathbb{T}^{n} \hookrightarrow ( \mathbb{R}^{2n}, \omega_{std} ) $ be a monotone Lagrangian embedding
   of the $n$-torus and let $ J $ be an arbitrary almost complex structure on $ \mathbb{R}^{2n}
   $, compatible with $ \omega_{std} $. Then for every point $ p \in L
   $ there exists a $ J $-holomorphic disc $ u: (D, \partial D)
   \rightarrow (M,L) $ whose boundary passes through $p$, i.e. $ p \in u( \partial D) $, and whose
   Maslov index $ \mu_{L}([u]) $ is $2$.

\end{thm}

Theorem~\ref{T:Audin conj} is a partial solution to a question of
Audin~\cite{Audin conj}, which states the same thing without the
monotonicity assumption. Previous results in this direction were
obtained by Biran and Cieliebak, Li, Oh, Polterovich, Viterbo
~\cite{Bi-Ci,Li,Oh-Spectral,A-L-P,P,V-1,V-2}.

 Our approach is based on Floer cohomology, in particular on the quantum product on Floer cohomology.
We will study the multiplicative behaviour of the spectral
sequence due to Oh, whose first page (i.e term $ E_1 $) is related
to the singular cohomology of a Lagrangian, and which converges to
its Floer cohomology.

 The idea of the proof of Theorem~\ref{T:Audin conj} is based on an idea originally raised by Seidel.
The proof of Theorem~\ref{T:Maslov 2} uses ideas
from~\cite{Bi-Co-2,Cornea-Lalonde}. The statement of
Theorem~\ref{T:Maslov 2} can be more directly proved (using Gromov
compactness theorem, but without any Floer theory) in the case of a
Clifford torus and of an exotic torus due to Chekanov
(see~\cite{Chekanov,E-P}), hence for every Lagrangian torus which is
Hamiltonianly isotopic to each one of them. However, the full
classification of monotone tori in $ \mathbb{R}^{2n} $ is still not
known.

  In this paper we use Floer cohomology with coefficients in a ring $ A = \mathbb{Z}_2 [T^{-1} , T] $,
and grade $ A $ by $ \deg T = N_{L} $ (see~\cite{Bi-Co-3,Oh-HF1}). We show

\begin{thm} \label{T:Proj space}

 Let $ L = \mathbb{RP}^{n}
   \hookrightarrow M $ be a monotone
   Lagrangian embedding of the real projective space into a tame symplectic manifold $ M $. Assume
   in addition that its minimal Maslov number is $ N_{L} \geqslant 3
   $. Then $ HF^*(L,L) \cong ( H(L; \mathbb{Z}_{2} ) \otimes A )^* $. In particular $ L
   $ is not Hamiltonianly displaceable. Moreover, for every
   Hamiltonian diffeomorphism $ f\colon M \rightarrow M $, for which $ f(L) $
   intersects $ L $ transversally, we have $ \sharp( L \cap f(L) ) \geqslant
   n+1 $.

\end{thm}

 After the first version of this paper was written, we received
from Fukaya, Oh, Ohta and Ono revised version of their
work~\cite{FOOO} in which results, similar to Theorem
~\ref{T:Audin conj}, were obtained.

 Section ~\ref{S:Proof of main results} is devoted for the proofs of
Theorems ~\ref{T:Audin conj}, ~\ref{T:Maslov 2}, ~\ref{T:Proj space}.
In Section ~\ref{S:Floer Spec} we recall the definitions of Floer
cohomology and the spectral sequence of Oh, and give the definition of
the quantum product on the Floer complex. Then, in the end of Section
~\ref{S:Floer Spec}, we state Theorems ~\ref{T:Leibnitz rule},
~\ref{T:E1}. Theorem ~\ref{T:Leibnitz rule} guarantees the
multiplicativity of the Floer cohomology and of the spectral sequence
of Oh. The statement of Theorem ~\ref{T:E1} is used in an essential
way in the proofs of Theorems ~\ref{T:Audin conj}, ~\ref{T:Maslov 2},
~\ref{T:Proj space}. Section ~\ref{S:Proofs of quantum product}
stands for the proofs of Theorems ~\ref{T:Leibnitz rule},
~\ref{T:E1}. Finally, Section ~\ref{S:basic notions} is devoted to
recall the basic notions from symplectic topology, that we use in the
present article.

\section{Proof of main results} \label{S:Proof of main results}
 In this section we provide proofs for Theorems ~\ref{T:Audin
 conj}, ~\ref{T:Maslov 2}, ~\ref{T:Proj space}. The tools that
we use in the proofs are the multiplicativity of the spectral
sequence of Oh, and special properties of this spectral sequence. For
detailed description of the Floer cohomology, the spectral sequence
of Oh, and the proof of existence of the multiplicative structure, we
refer the reader to Section ~\ref{S:Floer Spec}.

Denote by $ \{ E_{r}^{p,q} , d_{r} \} $ the spectral sequence of
Oh, with coefficients in the ring $ A = \mathbb{Z}_2 [T^{-1} , T]
$, where $ A $ is graded by $ \deg T = N_{L} $.
 The properties of it, that are essential for the proofs of Theorems
~\ref{T:Audin conj}, ~\ref{T:Maslov 2}, ~\ref{T:Proj space}, are (see
Section ~\ref{S:Floer Spec}) :

\begin{itemize}

  \item $E_{1}^{p,q}=H^{p+q-pN_{L}}(L, \mathbb{Z}_{2}) \otimes
   A^{pN_{L}}$ , $d_{1} = [ \partial_{1} ] \otimes T $, where $ [
   \partial_{1} ] : H^{p+q-pN_{L}}(L, \mathbb{Z}_{2}) \rightarrow
   H^{p+1+q-(p+1)N_{L}}(L, \mathbb{Z}_{2}) $ is induced from $
   \partial_{1} $ - the operator, that enters in the definition of
   the Floer differential (see Section ~\ref{S:Floer Spec}).
   The multiplication on $ H^{*} (L, \mathbb{Z}_{2}) $, induced
   from the multiplication on $E_{1}^{p,q}$, coincides with the
   standard cup product.

  \item More generally, for every $r\geq 1$, $E_r^{p,q}$ has the form $E_r^{p,q} =
   V_r^{p,q} \otimes A^{p N_L}$ with $d_r = \delta_r \otimes
   T^r$, where $V_r^{p,q}$ are vector spaces over $\mathbb{Z}_2$
   and $\delta_r$ are homomorphisms $\delta_r:V_r^{p,q} \to
   V_r^{p+r,q-r+1}$ defined for every $p,q$ and satisfy $\delta_r
   \circ \delta_r = 0$.  We have $$V_{r+1}^{p,q} = \frac{\ker
   (\delta_r: V_r^{p,q} \to
   V_r^{p+r,q-r+1})}{\textnormal{image\,} (\delta_r: V_r^{p-r,
   q+r-1} \to V_r^{p,q})}.$$

  \item $ \{ E_{r}^{p,q}, d_{r} \}$  converges to $ HF $, i.e
   $$ E_{\infty}^{p,q} \cong \frac{F^p HF^{ p+q }(L,L)}{F^{p+1} HF^{ p+q }(L,L)}
   ,$$ where $ F^p HF(L,L) $ is the filtration on $ HF(L,L) $,
   induced from the filtration $ F^p CF(L,L) $.

\end{itemize}

\begin{proof}[Proof of Theorem~\ref{T:Audin conj}]
   Consider a Lagrangian embedding $ L = \mathbb{T}^{n}
   \hookrightarrow ( \mathbb{R}^{2n}, \omega_{std} ) $, and assume
   by contradiction that $ N_{L} \geqslant3 $. The Floer cohomology $ HF(L,L) $
   is well-defined, and the spectral sequence of
   Oh, which computes it, becomes multiplicative.
   Look at the first stage $ E_{1} $ of the spectral sequence. We know that
   $V_1^{p,q} = H^{p+q-pN_L}(L;\mathbb{Z}_2)$, and the induced
   product on $ V_{1} $ is the classical cup-product and we have
   a differential $ \delta_{1} : V_{1} \rightarrow V_{1} $, which
   decreases the natural grading $ V_{1} $ by $ N_{L}-1$.
   The key observation now is that
   the entire cohomology ring $ H^{*}(L, \mathbb{Z}_{2}) $ is
   generated by the first cohomology $ H^{1}(L, \mathbb{Z}_{2}) $
   since $ L = \mathbb{T}^{n} $ . Therefore looking on the natural
   grading on the $ V_{1} $, the elements of degree $ 1 $ generate
   the whole $ V_{1} $. Now, because of $ \delta_{1} $ decreases
   the grading by $ N_{L} - 1 \geqslant 2 $, then for every $ a \in
   V_{1} $ of degree $ 1 $ we have $ \delta_{1}(a) = 0 $.  Then, since
   $ \delta_{1} $ satisfies the Leibnitz rule, the kernel $ \Ker(\delta_{1}) \subseteq
   V_{1} $ is a sub-ring, therefore $ \Ker(\delta_{1}) = V_{1} $,
   so we obtain $ \delta_{1} \equiv 0 $.
   Therefore $ V_{2} = H(V_{1},\delta_{1}) = V_{1} $, therefore
   $ V_{2} $ is isomorphic to $ V_{1} $ as a graded ring.
   Therefore we can apply the same argument for $ \delta_{2} : V_{2}
   \rightarrow V_{2} $, since $ \delta_{2} $ decreases the grading by $
   2N_{L} -1 \geqslant 2 $, so we will get that $ \delta_{2} \equiv 0 $,
   and so on, so at each stage we will get that $ \delta_{r} \equiv 0 $ by
   induction, therefore $ d_{1} = d_{2} = \ldots = 0 $, so
   as a conclusion we get that $ E_{1} = E_{2} =
   \ldots = E_{ \infty } = HF(L,L) $. But
   $ L $ is clearly Hamiltonianly displaceable in $ \mathbb{R}^{2n} $,
   hence $ HF(L,L)=0 $ and we obtain that $ E_{1} = 0 $. Contradiction.
\end{proof}

\begin{rem} \label{R:Audin generalized}
   The same arguments from the proof of Theorem~\ref{T:Audin conj} in
   fact prove the following more general result:

   Let $L \subset (M,\omega)$ be a monotone Lagrangian with $N_L\geq
   2$. Assume that:
     \begin{enumerate}
       \item $(M,\omega)$ is a subcritical Stein manifold. (See
        e.g.~\cite{Bi-NonIntersect} for the definition).
       \item $H^*(L;\mathbb{Z}_{2})$ is generated by
        $H^1(L;\mathbb{Z}_{2})$.
     \end{enumerate}
     Then $N_L=2$.
  \end{rem}

\begin{proof}[Proof of Theorem~\ref{T:Maslov 2}]

   First let us show this statement for generic
   $ J $ and then we will use a compactness argument for
   proving it for every $ J $, compatible with $ \omega_{std} $.
   Let us make a generic choice of $ J $ and a Morse function $ f:
   L \rightarrow L $, such that the point $ x \in L $ is the only
   point of maximum of $ f $, i.e the only point with index equal
   to $ n $ ( it is easy to show that such an $ f $ exists ).
   Then the Floer cohomology $ HF(L,L) $ is well-defined.
   By Theorem~\ref{T:Audin conj} we have
   that $ N_{L}=2 $. Look at the spectral sequence of Oh, which
   converges to $ HF(L,L) $. Let us show that the differential $ \delta_{1}
   : V_{1} \rightarrow V_{1} $ is non-zero. The argument is
   similar to the one from Theorem~\ref{T:Audin conj}.
   Indeed, if conversely, we have that $ \delta_{1} \equiv 0 $, then $ V_{2} =
   H(V_{1},\delta_{1}) = V_{1} $ as graded rings and $ \delta_{2}
   : V_{2} \rightarrow V_{2} $ is a differential which decreases
   the grading by $ 2N_{L}-1 = 3 $ and $ V_{2} = V_{1} $ is
   as a ring generated by it's elements of degree $ 1 $, therefore
   $ \delta_{2} = 0 $, so  $ V_{3} = H(V_{2},\delta_{2}) = V_{2} = V_{1} $
   as graded rings and so on. Therefore we will obtain a
   contradiction as in proof of Theorem~\ref{T:Audin conj} with
   the fact that $ HF(L,L) = 0 $. So we have proved that $ \delta_{1}
   : V_{1} \rightarrow V_{1} $ is non-zero, therefore the map $
   [ \partial_{1} ] : H^{*}(L, \mathbb{Z}_{2}) \rightarrow
   H^{*-1}(L, \mathbb{Z}_{2}) $ is non-zero. Let us show that
   this implies that $ [ \partial_{1} ] ( [p] ) \neq 0 $, where
   $ [p] \in H^{n}(L, \mathbb{Z}_{2}) $ is a generator. From the
   theorems above we have that $ [ \partial_{1} ] : H^{*}(L, \mathbb{Z}_{2}) \rightarrow
   H^{*-1}(L, \mathbb{Z}_{2}) $ is a differential and satisfies a
   Leibnitz rule. Therefore, since $ H^{1}(L, \mathbb{Z}_{2}) $
   generates the whole $ H^{*}(L, \mathbb{Z}_{2}) $, $ [ \partial_{1}
   ]$ restricted to $ H^{1}(L, \mathbb{Z}_{2}) $ is non-zero. Take
   some $ x_{1} \in H^{1}(L, \mathbb{Z}_{2}) $ such that $ [ \partial_{1}
   ](x_{1}) \neq 0 $, then because of $ [ \partial_{1}
   ](x_{1}) \in H^{0}(L,\mathbb{Z}_{2}) \cong \mathbb{Z}_{2} $ we
   have that $ [ \partial_{1} ] (x_{1}) = 1 \in
   H^{*}(L,\mathbb{Z}_{2}) $. Complete $ x_{1} $ to a basis $ x_{1},x_{2}, \ldots ,x_{n} $ of
   $ H^{1}(L, \mathbb{Z}_{2}) $ as a vector space over $ \mathbb{Z}_{2} $.
   Then it is easy to see that the product $ x_{1}x_{2} \ldots
   x_{n} = [p] \in H^{n}(L,\mathbb{Z}_{2}) $. Denote $ y =
   x_{2}x_{3} \ldots x_{n} $, then $ x_{1} ( [ \partial_{1} ]
   ([p])) = x_{1} ( [ \partial_{1} ](x_{1}y)) =
   x_{1} ( ([ \partial_{1} ](x_{1})) y + x_{1} ([ \partial_{1} ](y))) =
   x_{1}y + x_{1}^{2} ([ \partial_{1} ](y)) = x_{1}y = [p] \neq
   0 $, hence $ [ \partial_{1} ]([p]) \neq 0 $. Now if we go back
   to the definition of $ [ \partial_{1} ] $, we see that the
   moduli-space $ \mathcal{M}_{1}(p,q) $ is non-empty for some critical point $
   q $ of $ f $ with $ \ind_{f}(q) = n-1 $. Now, because of $ p $
   is a critical point to the top index, the only gradient
   trajectories which start with $ x $ are constant trajectories,
   therefore the boundary of the $ J $-holomorphic disc from the
   definition of $ \mathcal{M}_{1}(p,q) $ contains the point $ p
   $.

   We have proved the theorem for generic choice of $ J $. For the
   general $ J $, consider a sequence of generic $ J_{n} $ which
   converge to $ J $. Then for each $ J_{n} $ we have a $ J_{n} $
   -holomorphic disc $ u_{n} : (D, \partial D) \rightarrow ( \mathbb{R}^{2n},
   L) $ such that $ p \in u_{n}( \partial D ) $. By Gromov
   compactness theorem ( see~\cite{Gr} ) there is a subsequence of
   $ u_{n} $ which converges to a tree of discs with sphere
   bubblings, which are $ J $-holomorphic, however because of $ \mu(u_{n}) = N_{L} = 2 $ is
   minimal, in the limit we have only one disc and no bubblings of
   spheres occur, therefore we obtain a $ J $ -holomorphic disc
   which contains $ p $.

\end{proof}

\begin{proof}[Proof of Theorem~\ref{T:Proj space}]
   The idea is similar to the the one in the proof of
   Theorem~\ref{T:Audin conj}. Namely, as before we see that the Floer
   cohomology of $ L $ is well defined, and we can compute it via the
   spectral sequence of Oh, which is multiplicative.
   Look at the first stage $ E_{1} $ of the spectral sequence.
   We have that $ H^{1}(L, \mathbb{Z}_{2}) $
   generates the entire cohomology $ H^{*}(L, \mathbb{Z}_{2}) $,
   hence $ V_{1} $ is generated as a ring by it's elements of
   degree $ 1 $ and the differential $ \delta : V_{1} \rightarrow
   V_{1} $ decreases the grading by $ N_{L}-1 \geqslant 2 $. Then
   arguing as in the proof of Theorem~\ref{T:Audin conj}, we conclude
   that $ E_{1} = E_{2} = \ldots = E_{ \infty } = HF(L,L) = H^{*}(L;
   \mathbb{Z}_{2}) $. Therefore $ HF(L,L) \neq 0 $ , hence $ L $ is
   not displaceable. Moreover for every Hamiltonian diffeomorphism $
   f\colon M \rightarrow M $, for which $ f(L) $ intersects $ L $ transversally, we have $
   \sharp ( L \cap f(L) ) = \dim_{ \mathbb{Z}_{2} } CF(L;f(L))
   \geqslant \dim_{ \mathbb{Z}_{2} } H^{*}(L; \mathbb{Z}_{2}) = n+1 $.
\end{proof}

\section{Floer cohomology, the spectral sequence, and the quantum product} \label{S:Floer Spec}

\subsection{Floer cohomology and the spectral sequence of Oh}

We recall some basic facts of Floer theory, which will be stated
without proofs (for the proofs we refer the reader
to~\cite{Bi-Co-1,Bi-Co-2,Oh-Spectral}).  Let $(M, \omega)$ be a
tame symplectic manifold and let $L \subset (M, \omega)$ be a
monotone Lagrangian submanifold with minimal Maslov number $N_L
\geq 2$. Then one can define the Floer cohomology of the pair $
(L,L) $, which we will denote by $ HF(L,L) $. As mentioned before,
we will work here with coefficients in the ring $ A = \mathbb{Z}_2
[T,T^{-1}] $, as described in ~\cite{Bi-Co-3}. In fact, we will
work with an equivalent definition of $ HF(L,L) $ as described in
~\cite{Bi-Co-3,Oh-Quantum}, which uses holomorpic disks rather
than holomorphic strips. We briefly recall the construction now.

%Figure1

\begin{figure} [h]
 \begin{center}
  \includegraphics[scale=0.7]{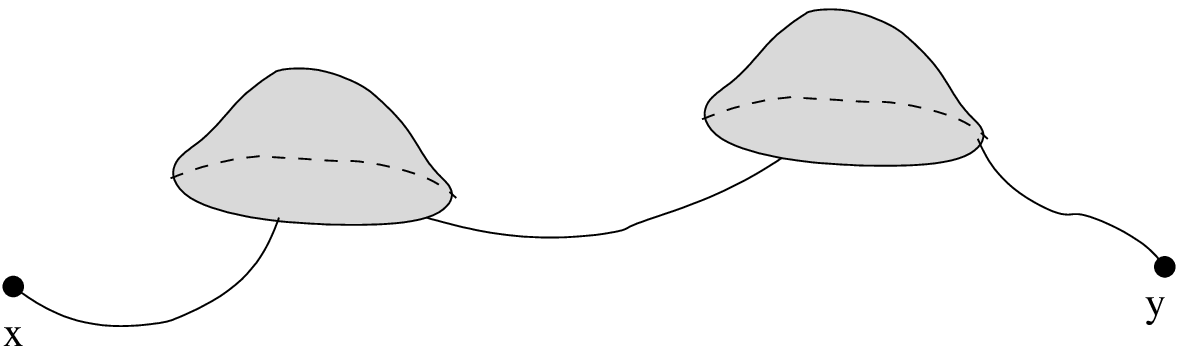}
 \end{center}
  \caption{}
   \label{F:fig1}
\end{figure}

We choose a generic pair of a Morse function $f\colon L \to
\mathbb{R}$ and a Riemannian metric on $ L $ and consider  a
generic  almost complex structure on $ M $. Denote by $ C_f^* $
the Morse complex associated to $ f $, graded by Morse indeces of
$f$. Let $ A = \mathbb{Z}_{2}[T,T^{-1}] $ be the algebra of
Laurent polynomials over $ \mathbb{Z}_{2} $ , where we take $ \deg
T = N_{L} $. Take the decomposition $ A = \bigoplus_{i \in
\mathbb{Z}} A^i $, where $ A^i $ is the subspace of homogenous
elements of degree $ i $.

 We define the Floer complex as $ CF(L,L) = C_{f} \otimes A $,
which has a natural grading coming from grading on $ C_f, A $.
More specifically, $$ CF^{l}(L,L) = \bigoplus_{k \in \mathbb{Z} }
C_f^{l-kN_{L}}(L,L) \otimes A^{kN_{L}} ,$$ for each $ l \in
\mathbb{Z} $. To define the Floer differential $ d_F : CF(L,L)
\rightarrow CF(L,L) $, we introduce auxiliary operators $
\partial_0,\partial_1,\dots,\partial_{\nu} : C_{f} \rightarrow
C_{f}$, where $\nu = \big[ \frac{\dim L + 1}{N_L} \big]$.

 For every pair of critical points
$ x,y \in C_{f} $ denote by $ \mathcal{M}_{k}(x,y) $ the
moduli-space of diagrams as in Figure~\ref{F:fig1}. This diagram
consists of pieces of gradient trajectories of $ -f $, joined by
somewhere injective(see ~\cite{L}) pseudo-holomorphic discs $ u_{1},u_{2}, \ldots
, u_{l} $ in $ M $, with boundaries on $ L $, such that the first
piece of gradient trajectory converges to $ x $, and the last
piece converges to $ y $, when time is reversed and such that sum
of Maslov indices of the discs is $ \mu(u_{1}) + \mu(u_{2}) +
\ldots + \mu(u_{l}) = kN_{L} $. Let us give a precise definition
of an element of $ \mathcal{M}_{k}(x,y) $. Consider a collection
of somewhere injective pseudo-holomorphic discs $ u_{j} \colon
(D^{2}, \partial D^{2}) \rightarrow (M,L) $, $ j=1,2, \ldots ,l$,
with boundaries on $ L $, and a collection of trajectories $$
\gamma_{0}: (-\infty , a_{0}] \rightarrow L , \gamma_{1}:
[a_{0},a_{1}] \rightarrow L , \gamma_{2}: [a_{1},a_{2}]
\rightarrow L , \ldots ,$$ $$ \gamma_{l}: [a_{l-1},a_{l}]
\rightarrow L, \gamma_{l+1}: [a_{l}, + \infty ) \rightarrow L ,$$
where $ -\infty < a_{0} < a_{1} < \ldots < a_{l} < +\infty $. Then
we demand that $ \gamma_{0}, \gamma_{1}, \ldots , \gamma_{l+1} $
are gradient trajectories for $ -f $ (with respect to our
Riemannian metric on $ L $), namely

$$ \dot{\gamma_{j}}(t) = - \nabla f ( \gamma_{j}(t) ) , \quad j= 0,1,2, \ldots ,
l+1 ,$$ with matching conditions :
$$ \lim_{t \rightarrow +\infty} \gamma_{0}(t) = x,\lim_{t \rightarrow -\infty} \gamma_{l}(t) = y, $$
$$ \gamma_{j}(a_{j}) = u_{j+1}(-1) , \quad j =0,1, \ldots , l-1 ,$$
$$ \gamma_{j}(a_{j-1}) = u_{j}(1) , \quad j =1,2, \ldots , l ,$$
and that the total Maslov of the collection $ u_{1},u_{2}, \ldots
,u_{l} $ of discs is $$ \mu( u_{1} ) + \mu( u_{2} ) + \ldots + \mu(
u_{l} ) = kN_{L} .$$ Then an element of $ \mathcal{M}_{k}(x,y) $ is
given by a collection $$ \{ (u_{1},u_{2}, \ldots ,u_{l}), (\gamma_{0},
\gamma_{1}, \ldots , \gamma_{l}) \} $$ as above, when we factorize it
by its group of inner automorphisms.

For generic choice of $ f $ and of almost complex structures (see
~\cite{Oh-HF1} for the details), $ \mathcal{M}_{k}(x,y) $ is a
manifold of dimension
$$ \dim \mathcal{M}_{k}(x,y) = \ind (y) - \ind(x) + kN_{L} - 1 $$
(see ~\cite{FOOO,Oh-Quantum}). If the dimension is $ 0 $, $
\mathcal{M}_{k}(x,y) $ is a manifold of dimension $ 0 $ and is
compact, so it is a finite collection of points. Denote in this
case $ n_{k}(x,y) := \sharp \mathcal{M}_{k}(x,y)(\mod 2) $. Then define
$$ \partial_{k}(x) = \sum_{ y \in C_{f} , ind(y) - ind(x) + kN_{L} - 1 =
  0 } n_{k}(x,y)y .$$
Note that $ \partial_0 $ is the usual Morse-cohomology boundary
operator. We see that $ d_{F} $ is not compatible with grading,
since each $
\partial_i $ acts like $\partial_i : C_{f}^* \to C_{f}^{*+1 - i
  N_L}(L)$.

 Let Floer differential $ d_F : CF^{*}(L,L) \rightarrow CF^{*+1}(L,L) $ by definition equals to
$ d_F = \partial_{0} \otimes 1 + \partial_{1} \otimes T + \ldots +
\partial_{\nu} \otimes T^{\nu} $. Then one can show that $
d_F $ is indeed a differential. It is well known that the homology of
the complex $ (CF(L,L),d_F) $ is canonically isomorphic to the Floer
cohomology of the pair $ (L,L) $ (see ~\cite{Bi-Co-1,Bi-Co-2} and the
references therein). Therefore we will write
$$ HF^*(L,L) = H^{*}(CF(L,L), d_F) .$$

 Consider now the following decreasing filtration on $ CF(L,L) $: $$
F^{p}CF(L,L) = \{ \sum x_{i} \otimes T^{n_{i}} \mid x_{i} \in
C_{f} , n_{i} \geqslant p \} .$$ It is obviously compatible with $
d_F $ (due to the monotonicity of $L$), so by a standard algebraic argument we obtain the spectral
sequence $ \{ E_{r}^{p,q}, d_{r} \}$. The following properties of
it have been proved in~\cite{Bi-NonIntersect}:

\begin{itemize}
  \item $ E_{0}^{p,q} = C_{f}^{p+q-pN_{L}} \otimes
   A^{pN_{L}}$ , $ d_{0} = [ \partial_{0} ] \otimes 1 $

  \item $E_{1}^{p,q}=H^{p+q-pN_{L}}(L, \mathbb{Z}_{2}) \otimes
   A^{pN_{L}}$ , $d_{1} = [ \partial_{1} ] \otimes T $, where $ [
   \partial_{1} ] : H^{p+q-pN_{L}}(L, \mathbb{Z}_{2}) \rightarrow
   H^{p+1+q-(p+1)N_{L}}(L, \mathbb{Z}_{2}) $ is induced from $
   \partial_{1} $.

  \item For every $r\geq 1$, $E_r^{p,q}$ has the form $E_r^{p,q} =
   V_r^{p,q} \otimes A^{p N_L}$ with $d_r = \delta_r \otimes
   T^r$, where $V_r^{p,q}$ are vector spaces over $\mathbb{Z}_2$
   and $\delta_r$ are homomorphisms $\delta_r:V_r^{p,q} \to
   V_r^{p+r,q-r+1}$ defined for every $p,q$ and satisfy $\delta_r
   \circ \delta_r = 0$.  Moreover $$V_{r+1}^{p,q} = \frac{\ker
   (\delta_r: V_r^{p,q} \to
   V_r^{p+r,q-r+1})}{\textnormal{image\,} (\delta_r: V_r^{p-r,
   q+r-1} \to V_r^{p,q})}.$$
   (For $r=0,1$ we have
   $V_0^{p,q}=C_f^{p+q-pN_L}$, $V_1^{p,q} =
   H^{p+q-pN_L}(L;\mathbb{Z}_2)$.)

  \item $ \{ E_{r}^{p,q}, d_{r} \}$ collapses at $ \nu + 1 $ step,
   namely $ d_{r} = 0 $ for every $ r \geqslant \nu +1 $, so $
   E_{r}^{p,q} = E_{\infty}^{p,q} $ for every $ r \geqslant \nu +1 $
   and the sequence converges to $ HF $, i.e
   $$ E_{\infty}^{p,q} \cong \frac{F^p HF^{ p+q }(L,L)}{F^{p+1} HF^{ p+q }(L,L)}
   ,$$ where $ F^p HF(L,L) $ is the filtration on $ HF(L,L) $,
   induced from the filtration $ F^p CF(L,L) $.

\end{itemize}

 Note that on each $ E_{r} $ and $ V_{r}
$ we have a natural grading coming from the grading on $ CF(L,L) $
and $ d_{r} : E_{r} \rightarrow E_{r} $ shifts this grading by $ 1
$. Therefore $ \delta_{r} : V_r \to V_r $ shifts the grading on $
V_r $ by $ 1-rN_{L} $, since the degree of $ T $ is $ N_{L} $. On
$ V_{1} $ this grading coincides with the natural grading on the
cohomology ring $ H^{*}(L, \mathbb{Z}_{2}) $.

\subsection{Quantum product on Floer cohomology}

 Consider generic Morse functions $ f,g,h : L
\rightarrow \mathbb{R} $ and denote by $$ ( CF_{f}, d_F^{f}),
(CF_{g}, d_F^{g}), (CF_{h}, d_F^{h}) $$ the corresponding Floer
complexes. Then we will be able to define a "quantum product" $
\star : CF_{f} \otimes CF_{g} \rightarrow
CF_{h} $, such that differentials will satisfy the
following analog of the Leibnitz rule: $$ d_F^{h} (a \star b) =
d_F^{f}(a) \star b + a \star d_F^{g}(b) .$$ Moreover, this product is
compatible with the filtrations on the $ CF $'s in the sense that $ * $
maps $ F^p CF_f \otimes F^{p'}CF_g $ to $ F^{p+p'}CF_h $.

%Figure2

\begin{figure} [h]
 \begin{center}
  \includegraphics[scale=0.7]{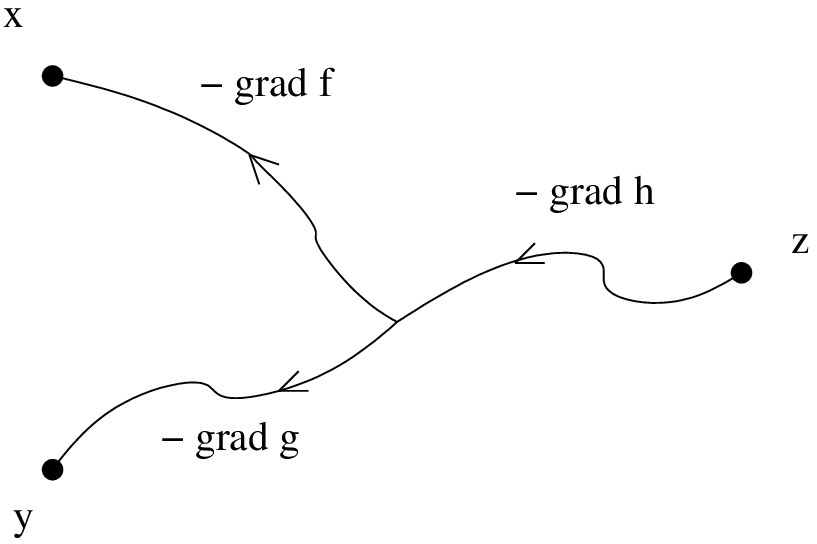}
 \end{center}
  \caption{}
   \label{F:fig2}
\end{figure}

Then automatically the spectral sequences become multiplicative, i.e we have products $
E_{r}^{p,q}(f) \otimes E_{r}^{p',q'}(g) \rightarrow
E_{r}^{p+p',q+q'}(h) $ at each stage of the spectral sequence,
which are induced from $ \star $, such that the differential $
d_{r} $ satisfies the Leibnitz rule, according to this product,
and such that the product at the $ r+1 $ stage comes from the
product at $ r $ stage. Note that in this case these products
induce products $ V_{r}^{p,q}(f) \otimes V_{r}^{p',q'}(g)
\rightarrow V_{r}^{p+p',q+q'}(h) $, the differential $
\delta_r:V_r^{p,q} \to V_r^{p+r,q-r+1} $ satisfies the Leibnitz
rule with respect to this product and that the product on the $
V_{r+1} $ is induced from the product on $ V_{r} $. Then the
crucial observation will be that on $ V_{1} $ the induced product
coincides with the usual cup-product on $ H^{*}(L; \mathbb{Z}_{2})
$, so the next products on the $ V_{r} $ are induced from this cup
product, therefore the quantum effects are lost for $ r \geqslant
1 $. Now let us describe the operation $ \star : F^{p}CF_{f}
\otimes F^{p'}CF_{g} \rightarrow F^{p+p'}CF_{h} $. To do this, let
us introduce operations $ m_{l} : C_{f} \otimes C_{g} \rightarrow
C_{h} $ (for a more general introduction of such an operations see
~\cite{Bets-Cohen,FOOO}).  The operation $ m_{0} $ is the usual
product on the Morse complexes of $ f,g,h $. Let us recall its
definition. For every triple of generators $ x \in C_{f},y \in
C_{g},z \in C_{h} $, denote by $ \mathcal{M}_{0}(x,y;z) $ the
moduli-space of diagrams as in Figure~\ref{F:fig2}. This diagram
is given by trajectories $ \gamma_{f} :[0, + \infty ) \rightarrow
L, \gamma_{g} :[0, + \infty ) \rightarrow L,\gamma_{h} :( -
\infty, 0] \rightarrow L,$ such that

$$ \dot{\gamma_{f}}(t) = - \nabla f (
\gamma_{f}(t)), \quad \dot{\gamma_{g}}(t) = - \nabla g (
\gamma_{g}(t)), \quad \dot{\gamma_{h}}(t) = - \nabla h (
\gamma_{h}(t)),$$
and
$$ \lim_{t \rightarrow + \infty} \gamma_{f}(t) = x,
\quad \lim_{t \rightarrow + \infty} \gamma_{g}(t) = y, \quad \lim_{t
  \rightarrow - \infty} \gamma_{h}(t) = z, $$
$$ \gamma_{f}(0)=\gamma_{g}(0)=\gamma_{h}(0).$$
Then by our assumption that $ f,g,h $ is a generic triple of Morse
functions, we get that $ \mathcal{M}_{0}(x,y;z) $ is a manifold,
and its dimension is given by $ \ind_{z}h - \ind_{x}f - \ind_{y}g
$. Moreover, when $\ind_{z}h - \ind_{x}f - \ind_{y}g = 0 $, $
\mathcal{M}_{0}(x,y;z) $ is a zero-dimensional compact manifold,
hence it consists of a finite number of points. For the case of $
\ind_{z}h - \ind_{x}f - \ind_{y}g = 0 $ we set $ n_{0}(x,y;z) :=
\sharp \mathcal{M}_{0}(x,y;z)(\mod 2) $. Now we define $$
m_{0}(x,y) := \sum_{ z \in C_{h}, \ind_{z}h -
  \ind_{x}f - \ind_{y}g = 0 } n_{0}(x,y;z)z .$$ This is the classical
cup product $ C_{f}^{i} \otimes C_{g}^{j} \rightarrow C_{h}^{i+j} $, and
so the classical Morse differential satisfies the Leibnitz rule and
induces the classical cup-product on cohomology $ H^{*}(L;
\mathbb{Z}_{2}) $. The further operations $ m_{l} $ for $ l \geqslant
1 $ will use quantum contributions.
%Figure3

\begin{figure} [h]
 \begin{center}
  \includegraphics[scale=0.7]{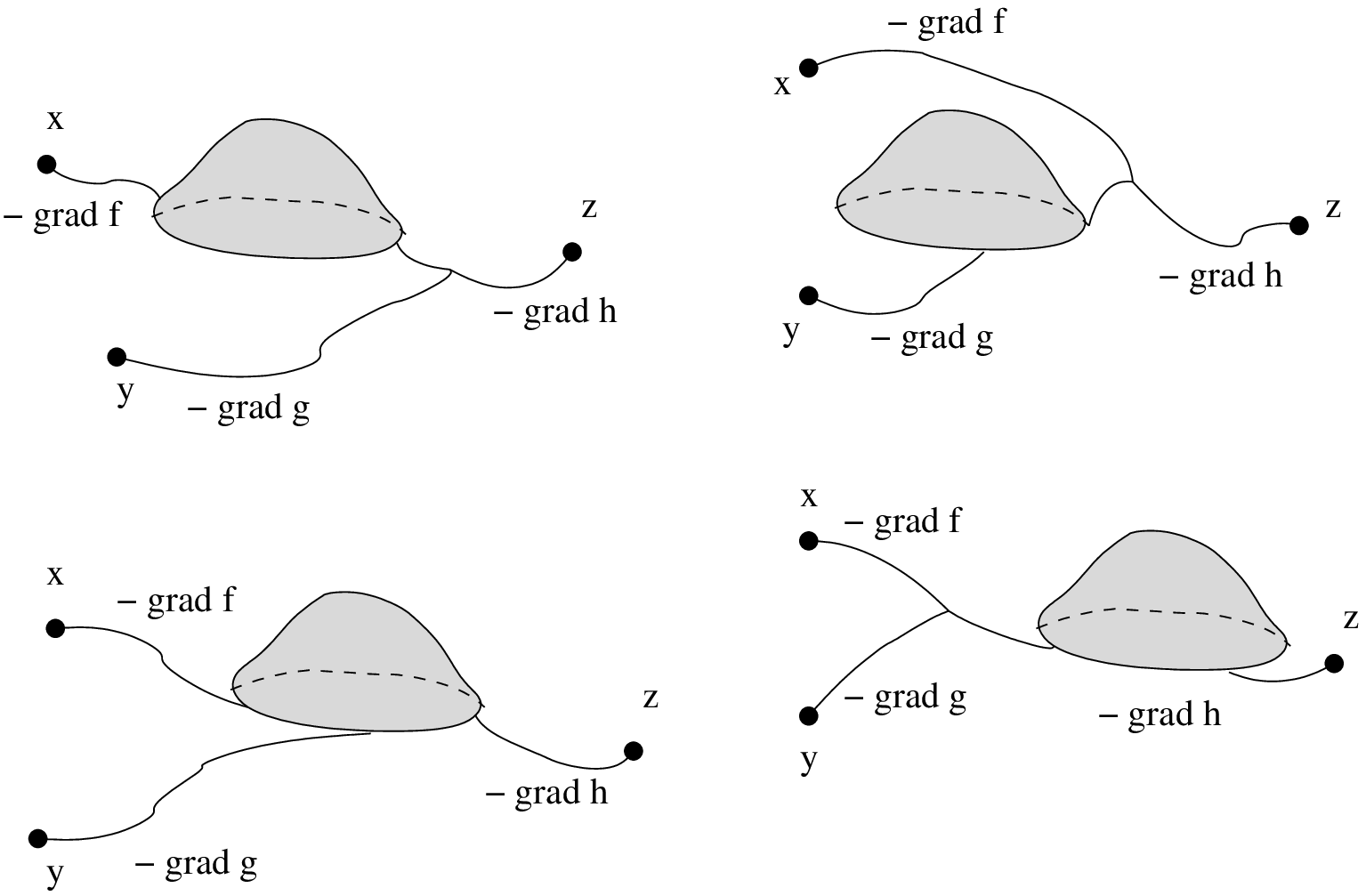}
 \end{center}
  \caption{}
   \label{F:fig3}
\end{figure}

Before defining the general $ m_{l} $, let us first describe $
m_{1} $. For a triple of critical points $ x \in C_{f},y \in
C_{g},z \in C_{h} $, we define the space $ \mathcal{M}_{1}(x,y;z)
$ to be the moduli-space of diagrams as in Figure~\ref{F:fig3},
where "black lines" are gradient trajectories of $ f,g,h $
respectively and the "discs" are pseudo-holomorphic somewhere
injective discs with boundaries on $ L $ and Maslov indices equal
to $ N_{L} $. Let us describe the first and the second diagram
from the Figure~\ref{F:fig3}. The first diagram is given by a
collection
$$ (u, \gamma_{f1},\gamma_{f2},\gamma_{g},\gamma_{h} ),$$ where
$ u: (D^{2},S^{1}) \rightarrow (M,L) $ is a pseudo-holomorphic
disc in $ M $ with boundary on $ L $ and $$ \gamma_{f1} :
[0,+\infty ) \rightarrow L , \gamma_{f2} : [0,a ] \rightarrow L ,
\gamma_{g} : [0,+\infty ) \rightarrow L , \gamma_{h} : ( -\infty,0
] \rightarrow L ,$$ where $ 0 < a \in \mathbb{R} $, such that
$$ \dot{\gamma}_{f1}(t)
= - \nabla f ( \gamma_{f1}(t)), \quad \dot{\gamma}_{f2}(t) = -
\nabla f ( \gamma_{f2}(t)), \quad \dot{\gamma_{g}}(t) = - \nabla g
( \gamma_{g}(t)),$$ $$ \dot{\gamma_{h}}(t) = - \nabla h (
\gamma_{h}(t)),$$

$$ \lim_{t \rightarrow + \infty} \gamma_{f1}(t) = x,
\quad \lim_{t \rightarrow + \infty} \gamma_{g}(t) = y, \quad \lim_{t
  \rightarrow - \infty} \gamma_{h}(t) = z, $$
$$ \gamma_{f2}(0)=\gamma_{g}(0)=\gamma_{h}(0), \quad
\gamma_{f1}(0)=u(-1), \quad \gamma_{f2}(a)=u(1). $$

The second diagram is given by a collection $ (u,
\gamma_{f},\gamma_{g},\gamma_{h} )$, where $ u: (D^{2},S^{1})
\rightarrow (M,L) $ is a somewhere injective pseudo-holomorphic
disc in $ M $, with boundary on $ L $ and $$ \gamma_{f} :
[0,+\infty ) \rightarrow L , \gamma_{g} : [0,+\infty ) \rightarrow
L , \gamma_{h} : ( -\infty,0 ] \rightarrow L ,$$ such that
$$ \dot{\gamma_{f}}(t) = - \nabla f (
\gamma_{f}(t)), \quad \dot{\gamma_{g}}(t) = - \nabla g (
\gamma_{g}(t)), \quad \dot{\gamma_{h}}(t) = - \nabla h (
\gamma_{h}(t)),$$

$$ \lim_{t \rightarrow + \infty} \gamma_{f}(t) = x,
\quad \lim_{t \rightarrow + \infty} \gamma_{g}(t) = y, \quad \lim_{t
  \rightarrow - \infty} \gamma_{h}(t) = z, $$
$$ \gamma_{f}(0)=u(1), \quad \gamma_{g}(0)=u(i), \quad \gamma_{h}(0)=u(-1). $$

\begin{figure} [h]
 \begin{center}
  \includegraphics[scale=0.7]{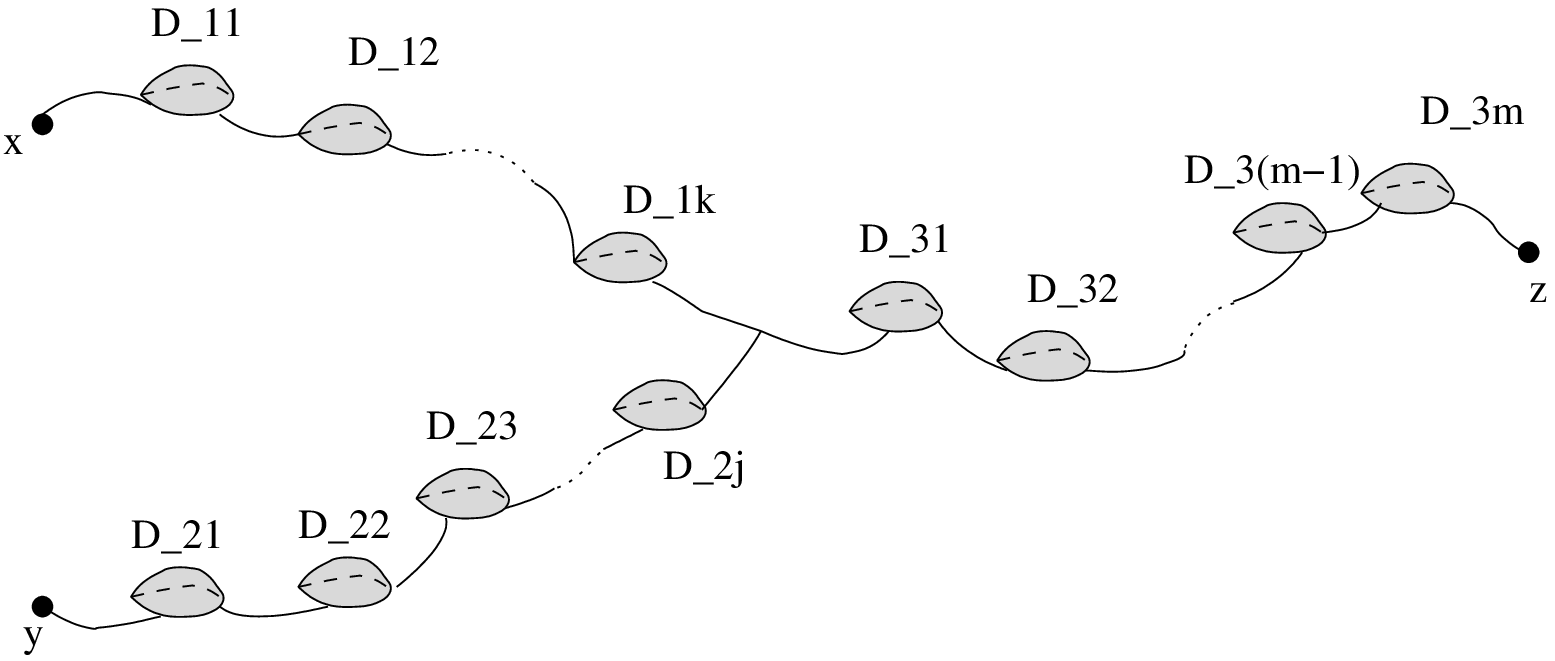}
 \end{center}
  \caption{}
   \label{F:fig4}
\end{figure}

\begin{figure} [h]
 \begin{center}
  \includegraphics[scale=0.7]{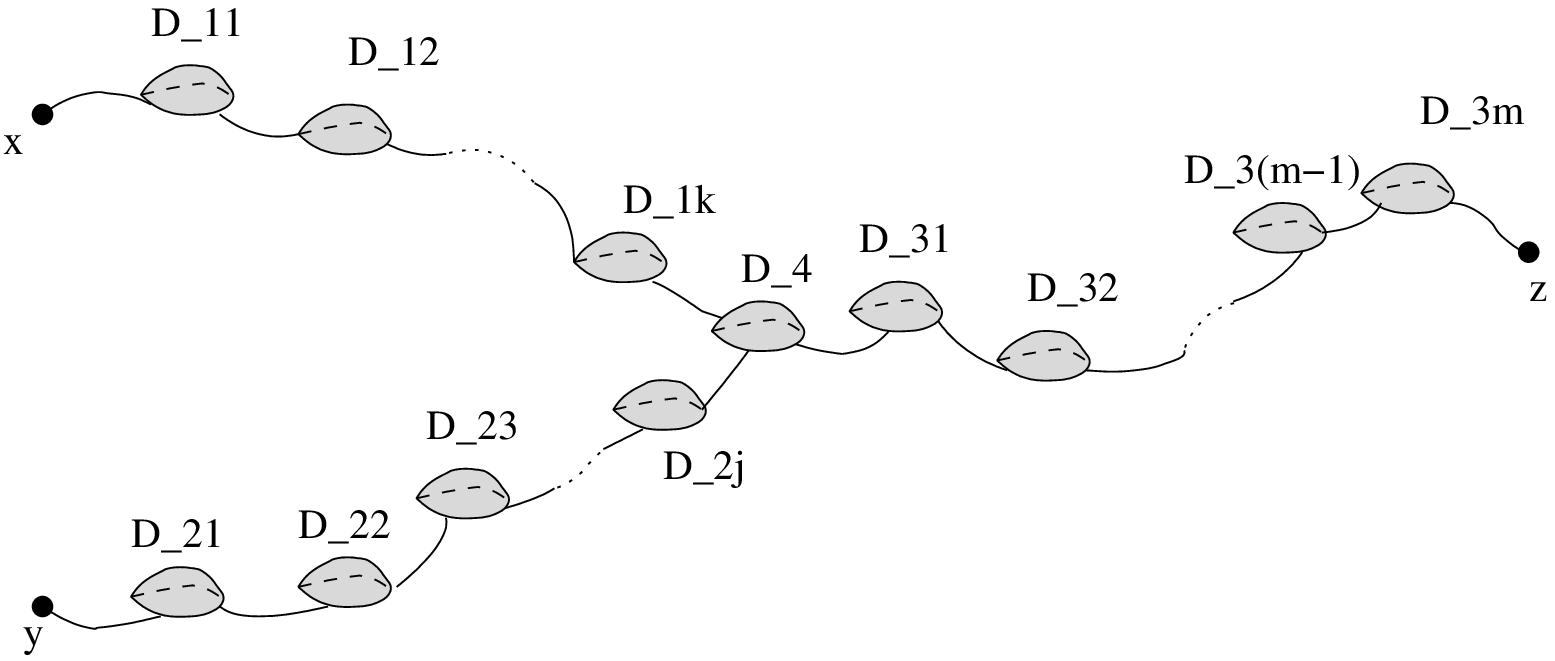}
 \end{center}
  \caption{}
   \label{F:fig5}
\end{figure}

The other diagrams from the picture have analogous definitions.  For
generic choices of $ f,g,h $ and of almost complex structures, $
\mathcal{M}_{1}(x,y;z) $ is a manifold of dimension $ \ind_{z}h -
\ind_{x}f - \ind_{y}g + N_{L} $. As before, if $ \dim
\mathcal{M}_{1}(x,y;z) = 0 $, it is a zero-dimensional compact
manifold, hence it is a finite collection of points and we set $
n_{1}(x,y;z) := \sharp \mathcal{M}_{1}(x,y;z) (\mod 2)$. Now we define
$$ m_{1}(x,y) := \sum_{ z \in C_{h},\ind_{z}h - \ind_{x}f - \ind_{y}g
  + N_{L} = 0 } n_{1}(x,y;z)z .$$ For defining $ m_{l} $ for general $
l $ we introduce a manifold $ \mathcal{M}_{l}(x,y;z) $ for $ x \in
C_{f},y \in C_{g},z \in C_{h} $. Its points are diagrams like in
Figure~\ref{F:fig3}, but instead they include several somewhere
injective discs with total Maslov index $ lN_{L} $. In general
these diagrams are of two types, as in Figures~\ref{F:fig4} and
~\ref{F:fig5}.

%Figures4,5

In Figure~\ref{F:fig4}, $ D_{11},D_{12}, \ldots , D_{1k} $ are
somewhere injective pseudo-holomorphic discs with boundaries on $
L $, which connect pieces of gradient trajectories of $ f $, $
D_{21},D_{22}, \ldots , D_{2j} $ - of $ g $, and $ D_{31},D_{32},
\ldots , D_{3m} $ - on $ h $, such that $ k,j,m \geqslant 0 $ and
total Maslov index
$$ \sum_{i,j} \mu(D_{ij}) = lN_{L} .$$
In Figure~\ref{F:fig5} in addition we have a disc $ D_{4} $ in the
middle, with gradient trajectory of $ h $ going into and gradient
trajectories of $ f,g $ going out of its boundary, and as before
the total Maslov index is
$$ \mu(D_{4}) + \sum_{i,j} \mu(D_{ij}) = lN_{L} .$$
We will denote by $ \mathcal{ \overline{M} }_{l}^{kjm} (x,y;z) $
the space of diagrams as in Figure~\ref{F:fig4} and by $ \mathcal{
\widetilde{M} }_{l}^{kjm} (x,y;z) $ the space of diagrams as in
Figure~\ref{F:fig5}. Denote also $ \mathcal{ \overline{M} }_{l}
(x,y;z) = \bigcup_{k,j,m} \mathcal{ \overline{M} }_{l}^{kjm}
(x,y;z) $, $ \mathcal{
  \widetilde{M} }_{l} (x,y;z) = \bigcup_{k,j,m} \mathcal{
  \widetilde{M} }_{l}^{kjm} (x,y;z) $. We have that $$
\mathcal{M}_{l}(x,y;z) = \mathcal{ \overline{M} }_{l}(x,y;z) \cup
\mathcal{ \widetilde{M} }_{l} (x,y;z)
$$ and $ \dim( \mathcal{ M }_{l}(x,y;z) ) = \ind_{z}h - \ind_{x}f -
\ind_{y}g + lN_{L} $.  Then, as before, in the case of $ \ind_{z}h -
\ind_{x}f - \ind_{y}g + lN_{L} = 0 $, we have that $ \mathcal{ M
}_{l}(x,y;z) $ is a finite collection of points and we set $
n_{l}(x,y;z) := \sharp \mathcal{ M }_{l}(x,y;z) $ in this case. Then,
as usual, we define
$$ m_{l}(x,y) :=  \sum_{ z \in C_{h}, \ind_{z}h - \ind_{x}f -
  \ind_{y}g + lN_{L} = 0 } n_{l}(x,y;z)z $$ for generators $ x \in
C_{f}, y \in C_{g} $.

Now we can define the quantum product $ \star : F^{p}CF_{f}
\otimes F^{p'}CF_{g} \rightarrow F^{p+p'}CF_{h} $ by $ x \star y =
m_{0}(x,y) \otimes 1 + m_{1}(x,y) \otimes T + m_{2}(x,y) \otimes
T^{2} + \ldots $ for $ x \in C_{f}, y \in C_{g} $, and then
naturally extend it to a map $ CF^{i}_{f} \otimes CF^{j}_{g}
\rightarrow CF^{i+j}_{h} $. Note that the filtrations on $ CF_{f},
CF_{g} ,CF_{h} $ are compatible with this map, i.e the image of $
F^{p}CF_{f} \otimes F^{p'}CF_{g}$ is in $ F^{p+p'}CF_{h} $. This
sum is finite, because $ \ind_{z}h - \ind_{x}f - \ind_{y}g +
lN_{L} = 0 $, hence $ l = ( -\ind_{z}h + \ind_{x}f + \ind_{y}g
)/N_{L} \leqslant 2n/N_{L} $. The main goal is now to prove the
following theorems.

\begin{thm} \label{T:Leibnitz rule} The differentials $
   d_F^{f},d_F^{g},d_F^{h} $ satisfy the
   Leibnitz rule with respect to the product $ \star $
   $$ d_F^{h} (a \star b) = d_F^{f}(a) \star b + a
   \star d_F^{g}(b), $$ for every $ a \in CF_{f},
   b \in CF_{g} $.
\end{thm}

\begin{thm} \label{T:E1}
   The product on $ V_{1} $, induced from $ \star $, coincides with the
   classical cup-product on $ H^{*}(L; \mathbb{Z}_{2}) $.
\end{thm}

\section{Existence and properties of the quantum product - proofs} \label{S:Proofs of quantum product}

\begin{proof}[Proof of Theorem~\ref{T:Leibnitz rule}]

   The main idea of the proof is similar to the analogous statement in
   classical Morse theory, for a standard Morse differential and a
   product on the Morse complex. In what follows we have to find a
   compactification of the manifolds $ \mathcal{ \overline{M} }_{l}
   (x,y;z) $ and $ \mathcal{ \widetilde{M} }_{l} (x,y;z) $. We are
   mostly interested in the components of the boundary of codimension
   $ 1 $. A point of the compactification of $ \mathcal{M}_{l}(x,y;z)
   $ is a diagram, consisting of several pseudo-holomorphic discs with
   boundaries on $ L $, spheres, critical points of $ f,g,h $ and
   pieces of gradient trajectories of $ f,g,h $ between them. Let us
   describe what can happen when we pass to the limit of a sequence of
   elements of $ \mathcal{M}_{l}(x,y;z) $. First, some of the gradient
   trajectories of $ f,g $ or $ h $ can "brake" and new critical
   points of $ f,g,h $ can appear in the diagram. Second, some of the
   pieces of the gradient trajectories can "shrink" to a point, such
   that we will obtain two touching discs or one disc containing a
   critical point. Also looking on $ \mathcal{ \overline{M} }_{l}
   (x,y;z) $, the piece of gradient trajectory containing the "middle
   point" can shrink to a point, so we will get a disc containing this
   "middle point". It can also happen that some of the discs split to
   a tree of discs and also some of the discs can bubble a sphere. The
   last thing that can happen is that for some pieces of trajectories
   which have endpoints on a disc, their endpoints on that disc
   converge one to another and become a single point. All these
   degenerations can happen simultaneously. However, when one looks on
   the codimension-$ 1 $ part of the boundary of $
   \mathcal{M}_{l}(x,y;z) $, only one of this degenerations can
   happen, moreover when we have the case of breaking of some
   trajectory, it can brake only at one point and this can happen only
   for one trajectory. Also, when we have "shrinking" of some
   trajectory, only one trajectory can shrink to a point.  Thirdly,
   when we have splitting of a disc to discs, then only one disc can
   split and only to two discs, and bubbling of spheres always has
   co-dimension $ \geqslant 2 $. Finally, when endpoints of some
   trajectories, lying on the boundary of some disc, become one point,
   we always have codimension $ \geqslant 2 $, except for the
   case when it is in $ \mathcal{ \widetilde{M} }_{l} (x,y;z) $, and
   two of the trajectories which end on the boundary of the "middle
   disc" in the limit have the same end on this boundary.

%Figures6,7

\begin{figure} [h]
 \begin{center}
  \includegraphics[scale=0.7]{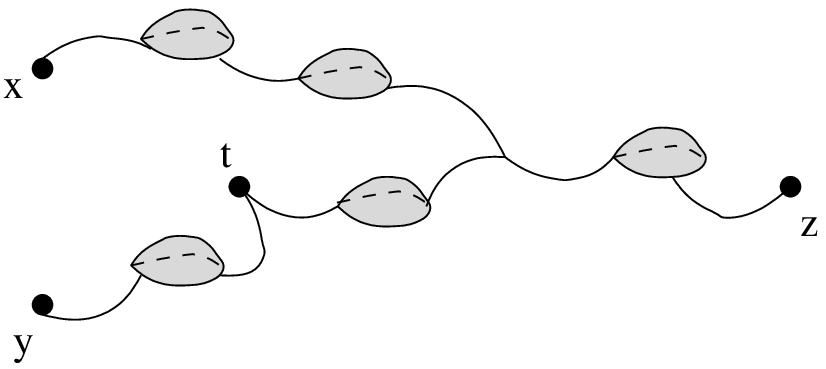}
 \end{center}
  \caption{}
   \label{F:fig6}
\end{figure}

\begin{figure} [h]
 \begin{center}
  \includegraphics[scale=0.7]{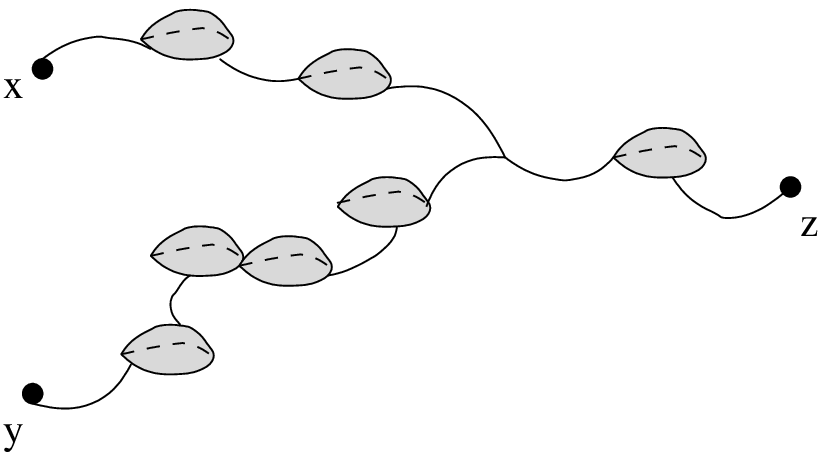}
 \end{center}
  \caption{}
   \label{F:fig7}
\end{figure}

    Therefore, when we are looking only at the codimension $ 1 $ part of $
   \mathcal{M}_{l}(x,y;z) $, only the following cases can happen:

   For $ \mathcal{ \overline{M} }_{l} (x,y;z) $ we can obtain:

   $a$) One trajectory "brakes" and we obtain a situation as in Figure
   ~\ref{F:fig6}, when a new critical point $ t $ appears ( this can happen with
   gradient trajectories of $ f,g,h $).

   $b$) Two neighboring discs in the chain become "touching" and the
   trajectory which joins them collapses into a
   point, as in Figure~\ref{F:fig7}.

   $c$) The last disc from the chain of discs of $ g $, for example,
   comes closer and closer to the "middle point", where gradient
   trajectories of $ f,g,h $ meet, until the trajectory which connects
   it to this "middle point" collapses to a point ( as in Figure~\ref{F:fig8} ).

   $d$) One of the discs splits to the union of two touching discs
   (Figure~\ref{F:fig9}).

%Figures8,9

\begin{figure} [h]
 \begin{center}
  \includegraphics[scale=0.7]{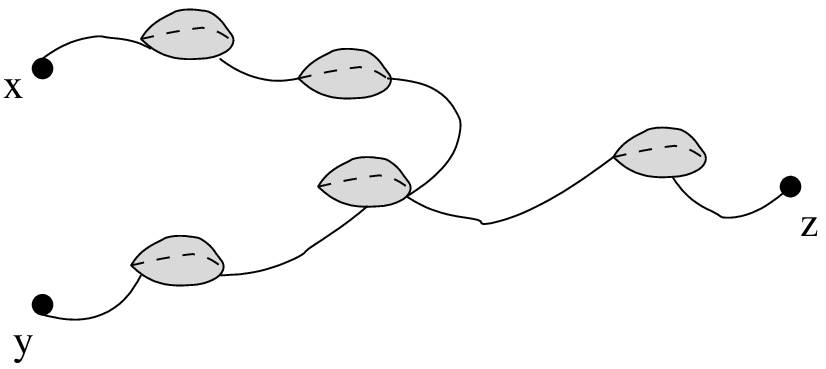}
 \end{center}
  \caption{}
   \label{F:fig8}
\end{figure}

\begin{figure} [h]
 \begin{center}
  \includegraphics[scale=0.7]{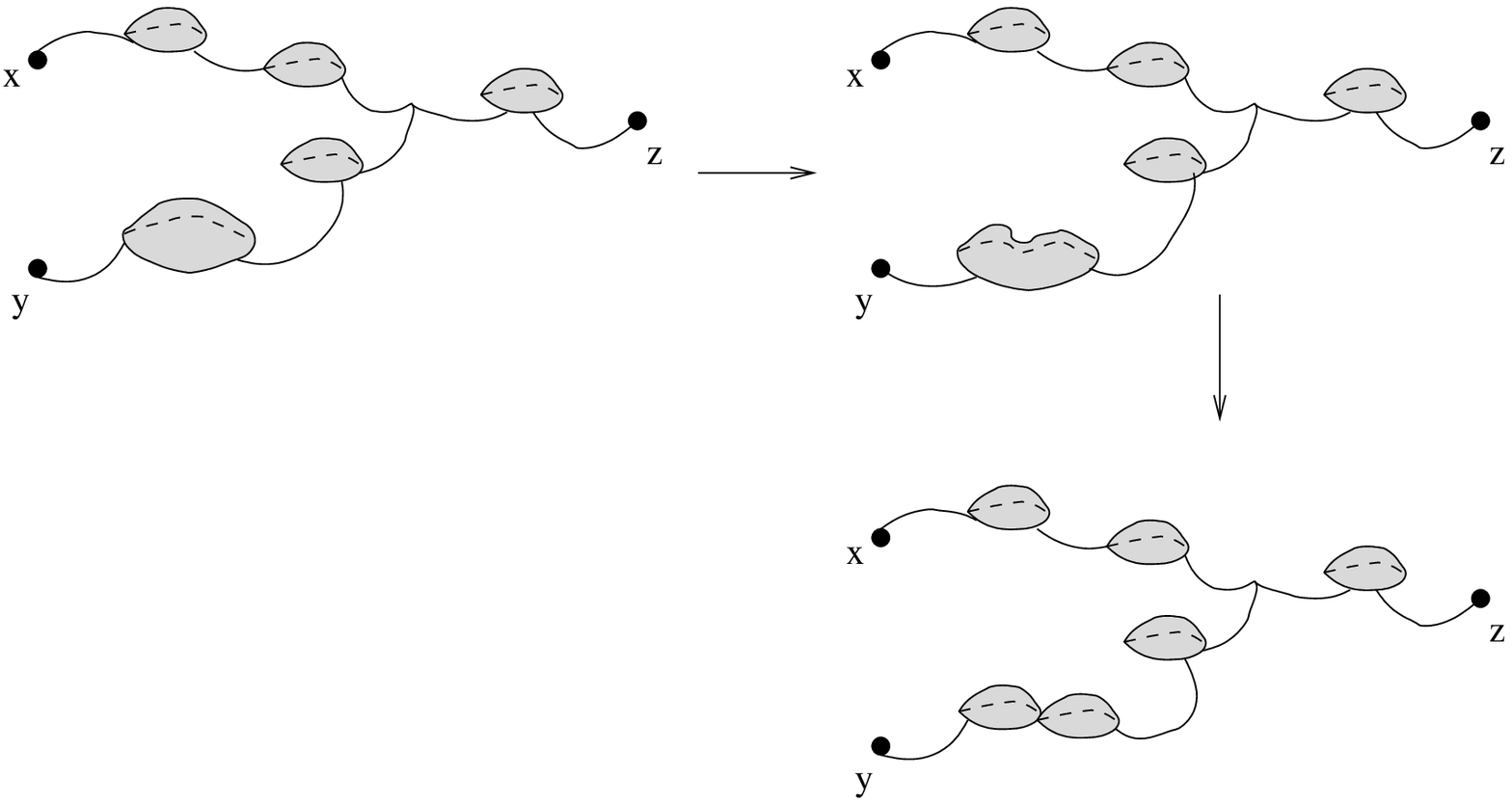}
 \end{center}
  \caption{}
   \label{F:fig9}
\end{figure}

   Denote by $$ \mathcal{M}^{a}_{l}(x,y;z),
   \mathcal{M}^{b}_{l}(x,y;z),
   \mathcal{M}^{c}_{l}(x,y;z),\mathcal{M}^{d}_{l}(x,y;z) $$ the
   manifolds of diagrams of types $ a,b,c,d $ respectively.

   Now addressing to the co-dimension $ 1 $ compactification of $
   \mathcal{ \widetilde{M} }_{l} (x,y;z) $, we see that the following
   cases are possible:

   $e$) One trajectory "brakes" and we obtain a situation as in Figure
   ~\ref{F:fig10}.

   $f$) Two neighboring discs in the chain come close and the
   trajectory which joins them collapses to a point, or the middle
   disc comes close to some neighboring disc as in Figure~\ref{F:fig11}.

   $g$) One of the discs splits to a union of two touching discs (and
   we again obtain a situation as in Figure~\ref{F:fig11}).

   $h$) Two trajectories touching the middle disc converge to
   trajectories which touch the middle disc in the same point.  We
   obtain again Figure~\ref{F:fig8}.

%Figures10,11

\begin{figure} [h]
 \begin{center}
  \includegraphics[scale=0.7]{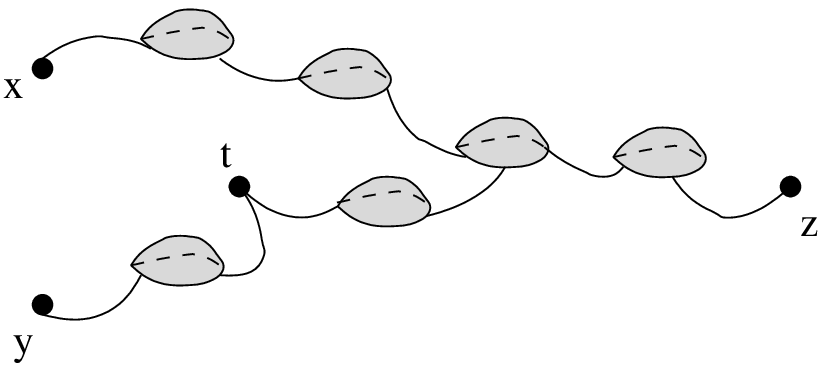}
 \end{center}
  \caption{}
   \label{F:fig10}
\end{figure}

\begin{figure} [h]
 \begin{center}
  \includegraphics[scale=0.7]{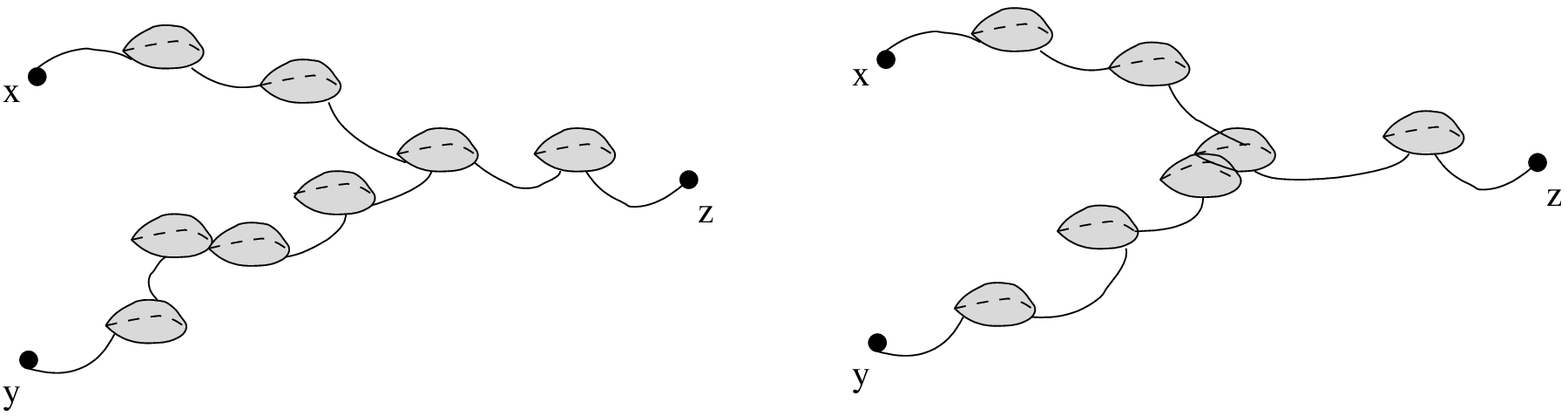}
 \end{center}
  \caption{}
   \label{F:fig11}
\end{figure}

   As before, we denote by $$ \mathcal{M}^{e}_{l}(x,y;z),
   \mathcal{M}^{f}_{l}(x,y;z),
   \mathcal{M}^{g}_{l}(x,y;z),\mathcal{M}^{h}_{l}(x,y;z)
    $$ the manifolds of the situations $ e,f,g,h $ respectively.
    Note that
    $$ (*)
      \begin{cases} %\label{eq:moduli}
         \mathcal{M}^{c}_{l}(x,y;z) = \mathcal{M}^{h}_{l}(x,y;z)  \\
         \mathcal{M}^{b}_{l}(x,y;z) = \mathcal{M}^{d}_{l}(x,y;z)  \\
         \mathcal{M}^{f}_{l}(x,y;z) = \mathcal{M}^{g}_{l}(x,y;z)
      \end{cases}
    $$
    Now let us see how this can be applied to prove the Leibnitz rule.
    Let us first write what it means. Taking $ x \in C_{f}, y \in
    C_{g} $, we have $$ d_F^{h}(x \star y) =
    d_F^{h} ( \sum_{ i \geqslant 0 } m_{i}(x,y) \otimes
    T^{i} ) = \sum_{ i,j \geqslant 0 } \partial_{j}( m_{i}(x,y) )
    \otimes T^{j+i} .$$ Similarly, $$ d_F^{f}(x) \star y =
    \sum_{ j,i \geqslant 0 } m_{i}(\partial_{j} x,y) \otimes T^{j+i}
    $$
    and $$
    x \star d_F^{g}(y) = \sum_{ j,i \geqslant 0 }
    m_{i}(x,\partial_{j} y) \otimes T^{j+i} .$$
    Therefore, we are left
    with proving that for every $ l $, $$
    \sum_{ i,j \geqslant 0, i+j
      = l } \partial_{j}( m_{i}(x,y) )= \sum_{ i,j \geqslant 0, i+j =
      l } m_{i}(\partial_{j} x,y) + \sum_{ i,j \geqslant 0, i+j = l }
    m_{i}(x,\partial_{j} y) .$$
    This means that we have to show that
    for every choice of generators $ x \in C_{f}, y \in C_{g}, z \in
    C_{h} $ with $ \ind(z) = \ind(x) + \ind(y) + lN_{L} +1 $, the
    total number of configurations in $ \mathcal{M}^{a}_{l}(x,y;z)
    \cup \mathcal{M}^{e}_{l}(x,y;z) $ is even. For this consider the
    space $ \mathcal{M}_{l}(x,y;z) $. It is a $ 1 $-dimensional
    manifold, therefore its boundary consists of an even number of
    points. On the other hand, from a gluing argument
    (see~\cite{Fukaya-Oh1,FOOO,McDuff-Jhol,Seidel}) it follows that
    this boundary is the union of

    $$ \mathcal{M}^{a}_{l}(x,y;z),
    \mathcal{M}^{b}_{l}(x,y;z),
    \mathcal{M}^{c}_{l}(x,y;z),\mathcal{M}^{d}_{l}(x,y;z), $$
    $$\mathcal{M}^{e}_{l}(x,y;z),
    \mathcal{M}^{f}_{l}(x,y;z),\mathcal{M}^{g}_{l}(x,y;z),
    \mathcal{M}^{h}_{l}(x,y;z)$$ (because in our case $ \dim
    \mathcal{M}_{l}(x,y;z)=1 $, so in a generic situation the part of
    the boundary of $ \mathcal{M}_{l}(x,y;z) $ of co-dimension bigger
    than $ 1 $ must be of dimension less than $ 0 $, so it is an empty
    set). Now, $ (*) $ shows that modulo $ 2 $, the total number of
    points on the boundary of $ \mathcal{M}_{l}(x,y;z) $ is equal to $
    \sharp \mathcal{M}^{a}_{l}(x,y;z) + \sharp
    \mathcal{M}^{f}_{l}(x,y;z) $ and is even. This proves
    Theorem~\ref{T:Leibnitz rule}.
 \end{proof}

\begin{rem} \label{R:Dimension formula}

     In several places we have applied the dimension formula $
     \dim\mathcal{M}(A,J) = n+\mu(A) $, where $ \mathcal{M}(A,J) $
     is a manifold of $ J $-holomorphic maps $ u \colon
     (D^{2}, \partial D^{2}) \rightarrow (M,L) $ with $ [u] = A
     \in \pi_{2}(M,L) $, in order to show that certain
     configurations of gradient trajectories and
     pseudo-holomorphic discs cannot appear for generic choice of
     $ J $ ( because of negative dimension ). However, this
     dimension formula is based on the transversality argument ( see~\cite{McDuff-Jhol} )
     which requires somewhere injectivity of the $ J $-holomorphic
     discs. To solve this problem we use work of Kwon, Oh and of
     Lazzarini ( see~\cite{K-O,L} ). More precisely, suppose we
     have such a configuration and some pseudo-holomorphic disc
     $ u : (D,\partial D) \rightarrow (M,L)  $ participated in it
     is not somewhere injective. Then we can decompose it to a
     union of almost everywhere injective
     discs $ u_{1},u_{2},\ldots, u_{k} $ with multiplicities $
     m_{1},m_{2},\ldots,m_{k} $ respectively, such that there
     exists $ m_{j} >1 $ ( see~\cite{K-O,L} for the precise details of this
     decomposition ). This decomposition preserves the relative homology class, namely
     $$ [u]=m_{1}[u_{1}] + m_{2}[u_{2}] + \ldots + m_{k}[u_{k}] \in  H_{2}(M,L) ,$$
     however, when we look on this configuration when we take all
     the discs $ u_{j} $ with multiplicity $ 1 $, then it's total
     area and hence a total Maslov class is strictly smaller than of the
     original configuration, therefore by a usual dimension count
     we obtain that such a configuration has negative dimension and
     hence cannot appear in a generic situation, therefore the original
     configuration also cannot appear. See ~\cite{Bi-Co-1} for more details on such arguments.
\end{rem}

\begin{proof}[Proof of Theorem~\ref{T:E1}]

     First look at the induced product on the level $ E_{0} $. Let us
     show that it coincides with the classical product on the Morse
     complex.  From the standard construction of the spectral sequence
     we have that
     $$ E_{0}^{p,q}(f) = F^{p} CF_{f}^{p+q} /  F^{p+1}
     CF_{f}^{p+q} \cong V_{0}^{p,q}(f) \otimes A^{p N_L}, $$
     $$E_{0}^{p',q'}(g) = F^{p'} CF_{g}^{p' + q'} /  F^{p'+1}
     CF_{g}^{p' + q'} \cong V_{0}^{p',q'}(g) \otimes A^{p' N_L},$$
     where $ V_{0}^{p,q}(f) = C_{f}^{p+q-pN_{L}} $,
     $ V_{0}^{p',q'}(g) = C_{g}^{p'+q'-p'N_{L}} $.
     Now take $ \alpha \in V_{0}^{p,q}(f)$, $ \beta \in
     V_{0}^{p',q'}(g) $. Now, we can take $ \overline{\alpha} :=
     \alpha \otimes T^{p} \in F^{p} CF_{f}^{p+q} $,
     $ \overline{\beta} := \beta \otimes T^{p'} \in F^{p'} CF_{g}^{p'+q'} $
     as a pre-images of $ \alpha \otimes T^{p} \in E_{0}^{p,q}(f)
     $, $ \beta \otimes T^{p'} \in E_{0}^{p',q'}(g) $ under
     natural projections $ F^{p} CF_{f}^{p+q} \rightarrow
     E_{0}^{p,q}(f)$, $ F^{p'} CF_{g}^{p'+q'} \rightarrow
     E_{0}^{p',q'}(g) $ respectively. Then, by definition of the product $
     \star $, $$ \overline{\alpha} \star \overline{\beta} = m_{0}( \alpha ,
     \beta) \otimes T^{p+p'} + m_{1}( \alpha ,
     \beta) \otimes T^{p+p'+1} + $$ $$ + m_{2}(
     \alpha , \beta) \otimes T^{p+p'+2} + \ldots
     \in F^{p+p'} CF_{h}^{q+q'},$$ and so the induced
     product of $ \alpha \otimes T^{p} \in E_{0}^{p,q}(f) $, $
     \beta \otimes T^{p'} \in E_{0}^{p',q'}(g) $
     is the image of $ \overline{\alpha} \star \overline{\beta}
     \in C_{h}^{p+p'+q+q'-pN_{L}-p'N_{L}} $ under the natural
     projection $ F^{p+p'} CF_{h}^{p+p'+q+q'} \rightarrow
     E_{0}^{p+p',q+q'}(h)$, which is $ m_{0}( \alpha ,
     \beta) \otimes T^{p+p'} $. Therefore the induced product of
     $ \alpha \in V_{0}^{p,q}(f)$, $ \beta \in
     V_{0}^{p',q'}(g) $ is $ m_{0}( \alpha , \beta)
     \in V_{0}^{p+p',q+q'}(h) $, which is the classical product in the Morse
     complex. Note also that the differential $ \delta_{0} :V_{0}^{p,q} \to
     V_{0}^{p,q+1} $ coincides with classical Morse differential, therefore
     the induced product on $ V_{1} = H( V_{0} , \delta_{0} ) $
     is the classical cup-product.
  \end{proof}

\section{Basic notions of symplectic topology in terms of Lagrangian Floer theory.} \label{S:basic notions}

 In this section we summarize some relevant notions from symplectic
topology used in the article.

\subsection{Tame symplectic manifold}
 A symplectic manifold $ (M, \omega) $ is called tame if there
exists an almost complex structure $ J $, such that the bilinear form
$ g_{ \omega , J }( \cdot,\cdot) = \omega( \cdot, J \cdot ) $ is a Riemmanian metric on $ M $,
and moreover the Riemmanian manifold $ (M , g_{ \omega , J } ) $ is geometrically bounded (i.e. its sectional curvature is bounded above and the injectivity radius is bounded below). See ~\cite{A-L-P,Gr} for more details and the relevance of this condition for the theory of pseudo-holomorphic curves.

\subsection{The Maslov class}

 The Maslov class is a homomorphism $\mu_{L}:\pi_{2}(M,L)
\rightarrow \mathbb{Z}$, associated to a Lagrangian submanifold $L
\subset (M,\omega)$ . To describe it, we start with the linear
case. Consider the space $\mathbb{R}^{2n} \cong \mathbb{C}^{n}$
with standard symplectic structure. Denote by $\mathcal{L}(n)$ the
set of all Lagrangian linear subspaces of $\mathbb{R}^{2n}$. The unitary
group $U(n)$ acts transitively on $\mathcal{L}(n)$ such that a stabilizer of
Lagrangian subspace $\mathbb{R}^{n}=\{ y_{1}=\ldots=y_{n}=0\}
\subset \mathbb{C}^{n}$ is the orthogonal group $O(n) \subset
U(n)$. So $\mathcal{L}(n)$ is homeomorphic to the quotient
$U(n)/O(n)$. On $U(n)/O(n)$ we have well-defined map $det^{2}:
U(n)/O(n) \rightarrow S^{1} \subset \mathbb{C}$, hence we obtain a
map $\mathcal{L}(n) \rightarrow S^{1}$. The corresponding
homomorphism $\mu: \pi_{1}( \mathcal{L}(n)) \rightarrow
\mathbb{Z}$ is called the Maslov index. It can be verified that it
is an isomorphism.

 Now, consider a symplectic manifold $(M,\omega)$ and a
Lagrangian submanifold $L \subset M$. Take a disc in $M$ with
boundary lying on $L$: $ u:(D,\partial D) \rightarrow (M,L)$. We
obtain the following diagram of vector bundles:
\begin{equation*}
 \begin{array}{clcr}
    u^{*}T(M) & \supset & u^{*}T(L) \\
    \downarrow &        & \downarrow \\
    D          &        & \partial D
 \end{array}
\end{equation*}
 Over each point on the disc $D$ we have a symplectic linear space
and for every point on the boundary $\partial D$ we have a
Lagrangian linear subspace of the corresponding linear symplectic
space. Now, we can symplectically trivialize the bundle $u^{*}T(M)
\rightarrow D$, and as a result we will get a loop of Lagrangians
in $\mathbb{R}^{2n}$, $\gamma: S^{1} \rightarrow \mathcal{L}(n)$.
Applying to this loop the Maslov index, we get the Maslov class
evaluated on $D$, namely, we define $\mu_{L}(D)=\mu(\gamma)$. It
can be shown that this definition does not depend on the
trivialization, and actually depends only on $[D] \in
\pi_{2}(M,L)$. Given a Lagrangian submanifold $L \subset
(M,\omega)$ we denote by $N_{L}$ the positive generator of the
image $\mu_{L}(\pi_{2}(M,L)) \subset \mathbb{Z}$. We shall refer
to $N_{L}$ as the minimal Maslov number of $L$.

\subsection{Symplectic area class}

 This is a homomorphism $\omega : \pi_{2}(M,L) \rightarrow \mathbb{R}$
which computes the symplectic area of a disc: if we have a
representative $\alpha=[\varphi : (D , \partial D ) \rightarrow
(M,L)] \in \pi_{2}(M,L)$ , we define $\omega(\alpha):= \int_{D}(
\varphi^{*}\omega) $. Also here it can be shown that the
symplectic area class $\omega(\alpha)$ depends only on  $\alpha
\in \pi_{2}(M,L)$.

\subsection{Monotone symplectic manifolds}

 Let $(M,\omega)$ be a symplectic manifold. Denote by
$c_{1} \in H^{2}(M;\mathbb{R})$ the first Chern class of the
tangent bundle $T(M)$ viewed as a complex vector bundle (where the
complex structure on $T(M)$ is taken to be any almost complex
structure tamed by $\omega $). We say that $(M,\omega)$ is a
monotone symplectic manifold if there exists a positive real
number $\lambda>0$ such that $[\omega]=\lambda c_{1}$ .
 Given a Lagrangian submanifold $L \subset (M,\omega )$, we
have two homomorphisms: $\mu_{L}:\pi_{2}(M,L) \rightarrow
\mathbb{Z}$, and $\omega:\pi_{2}(M,L) \rightarrow \mathbb{R}$. We
say that $L$ is monotone, if these two homomorphisms $\mu_{L}
, \omega$ are proportional by some positive constant, that is,
there exists a constant $\lambda \in \mathbb{R}, \lambda > 0$ ,
such that for every $\alpha \in \pi_{2}(M,L)$ we have
$\mu_{L}(\alpha)=\lambda \omega(\alpha) $.

\subsubsection*{Acknowledgments}
I would like to thank my supervisor Paul Biran for his help and
attention he gave to me. I am grateful to Felix Schlenk for his
comments and for helping me to improve the quality of the
exposition. Also I am grateful to Leonid Polterovich, Alex Ivri
and Laurent Lazzarini for useful comments.

\end{document}